\documentclass[%
   , colorlinks 
]{mpi2015-cscpreprint}

\usepackage[american]{babel}
\usepackage[T1]{fontenc}
\usepackage[utf8]{inputenc}

\usepackage{graphicx}

\usepackage{amsmath}
\usepackage{amssymb}
\usepackage{amsthm}
\usepackage{cleveref}

\AtBeginDocument{
  \newtheorem{theorem}{Theorem}[section]
  
  \theoremstyle{remark}
  \newtheorem{remark}[theorem]{Remark}
}

\begin{document}


\title{Partial Observation of Linear Systems with
the Mori-Zwanzig Formalism}

\author[$*,\ddag$]{Fan Wang}
\affil[$\ddag$]{Department of Mathematics, Otto-von-Guericke-University, Magdeburg, Germany.\authorcr%
\email{fwang@mpi-magdeburg.mpg.de}}
\author[$\ast$]{Peter Benner}
\affil[$\ast$]{Computational Methods in Systems and Control Theory, Max Planck Institute for Dynamics of Complex Technical Systems, Magdeburg, Germany.\authorcr%
    \email{benner@mpi-magdeburg.mpg.de}}
\author[$\dagger$]{Jan Heiland}
\affil[$\dagger$]{Department of Mathematics and Natural Sciences, TU Ilmenau, Ilmenau, Germany.\authorcr%
  \email{jan.heiland@tu-ilmenau.de}}

\shorttitle{Partial Observation with MZ Formalism}
\shortauthor{Fan Wang, Peter Benner, Jan Heiland}
\shortdate{}



\abstract{%
The Mori-Zwanzig formalism provides a systematic framework for deriving reduced-order model of dynamical systems when only part of the state is observed, but its practical use is often limited by the complexity of the resulting computations. This paper develops an explicit formulation of the Mori-Zwanzig equation for linear time-invariant systems under partially observed observables. By expressing the dynamics in terms of observables, the Koopman generator, and projections onto resolved and unresolved components, we derive closed-form representations of the Markovian, noise, and memory contributions that arise in the Mori-Zwanzig identity.

For the linear setting, the resulting formulas recover the reduced dynamics obtained from the variation-of-constants formula while retaining the operator-based structure of the Mori-Zwanzig approach. This makes the derivation a transparent reference case for reduced-order modelling with memory and clarifies how unresolved variables influence the observed dynamics through history-dependent terms. The analysis also identifies the ingredients needed for extensions to nonlinear systems and more general projections. Analytical and numerical examples involving the harmonic oscillator and wave equations illustrate the construction and demonstrate how the formalism can be used to obtain interpretable reduced-order models for partially observed systems.}

\novelty{\par
\begin{itemize}
\setlength{\itemsep}{0.2em}
\setlength{\parskip}{0pt}
\setlength{\parsep}{0pt}
\item Derivation of an explicit Mori-Zwanzig formulation for partially
observed linear systems in terms of resolved and unresolved
state components.
\item Closed-form identification of the Markovian, noise, and memory
contributions, together with the associated memory kernel, for the resolved
dynamics.
\item Explicitly derive Mori-Zwanzig equation for the harmonic oscillator as an
illustrative example, showing how Markovian, noise, and 
memory terms arise.
\end{itemize}}

\maketitle


\section{Introduction}

Model Order Reduction (MOR) techniques, such as Proper Orthogonal
Decomposition (POD), Dynamic Mode Decomposition (DMD), and deep autoencoders, have enabled the efficient representation of complex dynamical systems. By identifying low-dimensional subspaces that capture the dominant features of a system, these methods significantly reduce computational costs while keeping the relative good accuracy of high-fidelity models. However, standard MOR methods face challenges in scenarios of partial observation, where only a subset of the system variables can be directly measured or resolved.

Two primary approaches have emerged to address the challenge of partial observation. The first relies on Takens' Embedding Theorem \cite{Tak81a}, which establishes that the full dynamics of a system can be topological identically reconstructed using time-delayed measurements of a single observed variable. This theorem forms the foundation of state-space reconstruction techniques widely used in nonlinear time series analysis. Complementing this theoretical framework, reservoir computing \cite{Jae01} has demonstrated remarkable success as a data-driven approach. By mapping inputs into a high-dimensional, nonlinear dynamical reservoir, these methods effectively capture complex temporal dependencies without requiring explicit knowledge of the unobserved variables.

While both Takens' embedding and reservoir computing can preserve the global dynamic properties of the original system, such as Lyapunov exponents, they do not necessarily guarantee the accurate reconstruction of the exact trajectories of the state variables. To address this limitation, we turn to the Mori-Zwanzig (MZ) formalism.

The MZ formalism provides a mathematically rigorous framework for dimensionality reduction. By projecting the full system dynamics onto a lower-dimensional subspace of ``resolved'' variables, the MZ decomposition yields the Generalized Langevin Equation. This exact reformulation consists of three distinct terms: a Markovian term (instantaneous dynamics), a memory kernel (non-Markovian history dependence), and an orthogonal dynamics term (often treated as noise). Unlike empirical approximations, the MZ formalism theoretically ensures the exact evolution of the resolved variables \cite{Vrugt2020}.

Recent research has demonstrated the practical power of the MZ framework in complex systems, particularly in turbulence modeling. Parish and Duraisamy \cite{Parish2017} utilized the formalism to derive non-Markovian closure models for Large Eddy Simulation (LES), demonstrating that explicit memory history improves the prediction of subgrid-scale energy transfer. Building on this, they developed a dynamic subgrid-scale model \cite{Parish2017Dynamic} that allows the memory length to be determined adaptively, and proposed a unified framework combining MZ with the Variational Multiscale Method to address scale separation issues more robustly \cite{Parish2017Unified}. Furthermore, Gouasmi et al. \cite{Gouasmi2017} established a priori estimation techniques to quantify these memory effects, verifying that the MZ formalism offers a viable path for physics-based closures in computational fluid dynamics.

Parallel to these physics-based advances, the intersection of MZ with
data-driven methods has gained significant traction. Lin et al.\cite{Lin2021}
formalized the connection between MZ and the Koopman operator, proposing a ``Deep
Mori-Zwanzig'' framework. By viewing the MZ formalism as a generalization of the
Koopman learning framework, they successfully learned memory kernels directly
from data. More recently, Lin et al. have expanded this by introducing regression-based projection operators \cite{Lin2023}, which allow for the use of nonlinear ansatzes (such as neural networks) to approximate the optimal Zwanzig projection. This line of inquiry has led to the development of the Mori-Zwanzig Mode Decomposition (MZMD) \cite{Woodward2024}, a novel technique that outperforms standard DMD in transient flows by explicitly accounting for memory effects in the modal decomposition. Despite its theoretical elegance and recent successes, the practical implementation of the MZ formalism remains challenging. The primary bottleneck is the computation of the memory kernel and the orthogonal dynamics, which typically requires solving the full high-dimensional system.

In this paper, we explore the derivation of the Mori-Zwanzig formalism
specifically for linear systems. While the ultimate goal of the formalism is to
address complex nonlinear dynamics by only accessing partially observed data. In
nonlinear systems, the memory kernel and orthogonal dynamics are generally
analytically intractable, often requiring simplification or numerical
approximations. In contrast, linear systems allow for the exact, closed-form
derivation of these terms. This provides algorithms and shows the precise
mechanisms by which projection operators transform high-dimensional couplings
into non-Markovian memory and noise. By mastering the exact structure of the MZ
decomposition in this tractable setting, we establish a rigorous foundation for
high dimensional linear Ordinary Differential Equations (ODEs) and Partial
Differential Equations (PDEs). Additionally, with the explicit formulation of
the relevant terms at hand, we provide a basis for future developments in more
complicated applications. The structure of this paper is organized as follows:
Section \ref{sec:MZ_formalism} reviews the theoretical foundation of MZ
formalism; Section \ref{sec:linear_MZ} derives the exact MZ formalism for linear
systems with projection operators that single out the resolved variables;
Section \ref{sec:Numerical_example} demonstrates numerical experiments on ODEs
and spatially discretized PDEs; and Section \ref{sec:conclusion} discusses implications and future directions.

\section{The Mori-Zwanzig Formalism} \label{sec:MZ_formalism}
In this section, we outline the mathematical foundations of the Mori-Zwanzig formalism and present a systematic derivation for projecting high-dimensional dynamical systems onto a low-dimensional subspace spanned by the observed data. We begin by defining the dynamical system and the associated operator-theoretic framework, then introduce the projection operators, and finally derive the Mori-Zwanzig Equation (MZE) via the Dyson identity, that is, the generalized Langevin equation that captures both memory effects and residual noise.

\paragraph{Dynamical Systems and Operators.}

Consider a general autonomous dynamical system defined on a smooth manifold $\mathcal{X} \subseteq \mathbb{R}^N$. The evolution of the state vector $\boldsymbol{x}(t) \in \mathcal{X}$ is governed by the ordinary differential equation:
\begin{equation*}
    \frac{d \boldsymbol{x}(t)}{d t} = \boldsymbol{R}(\boldsymbol{x}(t)), \quad \boldsymbol{x}(0) = \boldsymbol{x}_{0},
\end{equation*}
where $\boldsymbol{R}: \mathcal{X} \to \mathcal{X}$ is a vector field representing the system dynamics.

We are often interested in the evolution of specific quantities of interest rather than the full state $\boldsymbol{x}$. Let $\psi: \mathcal{X} \to \mathbb{R}$ be a real-valued observable and let $\boldsymbol{\xi}\in\mathcal{X}$ denote a generic point in state space. We consider a vector of observables $\boldsymbol{\psi}(\boldsymbol{\xi}) = [\psi_{1}(\boldsymbol{\xi}), \dots, \psi_{m}(\boldsymbol{\xi})]^\top$. The space of observables is defined as the set of smooth functions $\mathcal{S} = C^\infty(\mathcal{X})$.

Notably, the time evolution of an observable $\psi(\boldsymbol{x}(t))$ along the trajectory of the system can be described using the Koopman operator formalism. By the chain rule, the time derivative is given by

\begin{equation*}
    \frac{d}{dt} \boldsymbol{\psi}(\boldsymbol{x}(t))  = \nabla_{\boldsymbol{x}} \boldsymbol{\psi}(\boldsymbol{x}(t)) \cdot \dot{\boldsymbol{x}}(t) = \nabla_{\boldsymbol{x}} \boldsymbol{\psi}(\boldsymbol{x}(t)) \cdot \boldsymbol{R}(\boldsymbol{x}(t)),
\end{equation*}

\noindent which defines the so-called Liouville operator $\mathcal{L}: \mathcal{S} \to \mathcal{S}$ as the generator of the dynamics
\begin{equation*}
    (\mathcal{L}\psi)(\boldsymbol{\xi}) = \sum_{k=1}^N R_k(\boldsymbol{\xi}) \frac{\partial \psi}{\partial \xi_{k}}(\boldsymbol{\xi}) = \boldsymbol{R}(\boldsymbol{\xi}) \cdot \nabla_{\boldsymbol{\xi}}\psi(\boldsymbol{\xi}), \qquad \mathcal{L}\boldsymbol{\psi} = \begin{bmatrix} \mathcal{L}\psi_1 \\ \vdots \\ \mathcal{L}\psi_m \end{bmatrix}.
\end{equation*}
Consequently, the evolution of the observable vector satisfies the linear partial differential equation $\frac{d}{dt}\boldsymbol{\psi} = \mathcal{L}\boldsymbol{\psi}$. The formal solution is given by the semigroup generated by $\mathcal{L}$, denoted by $e^{t\mathcal{L}}$ and also known as the Koopman operator:
\begin{equation*}
    \boldsymbol{\psi}(\boldsymbol{x}(t;\boldsymbol{x}_{0}))
    = (e^{t \mathcal{L}} \boldsymbol{\psi})(\boldsymbol{x}_{0}).
\end{equation*}

\paragraph{Derivation of the Mori-Zwanzig Equation.}

The Mori-Zwanzig formalism seeks to partition the dynamics of $\boldsymbol{\psi}$ into a resolved part (projected onto a specific subspace) and an unresolved (orthogonal) part. We introduce a projection operator $\mathcal{P}$ acting on the space of observables, and its complement $\mathcal{Q} = \mathcal{I} - \mathcal{P}$.

To separate the dynamics, we employ the Dyson identity \cite{Evans2007-po}, following the standard Mori-Zwanzig derivation \cite{Vrugt2020}. For any two linear operators $A$ and $B$, the exponential $e^{t(A+B)}$ satisfies:
\begin{equation*}
    e^{t(A+B)} = e^{tB} + \int_0^t e^{(t-s)(A+B)} A e^{sB} ds.
\end{equation*}
We decompose the Liouville operator as $\mathcal{L} = \mathcal{PL} + \mathcal{QL}$. Applying the Dyson identity with $A = \mathcal{PL}$ and $B = \mathcal{QL}$, we obtain the operator identity:
\begin{equation}
    e^{t\mathcal{L}} = e^{t\mathcal{QL}} + \int_{0}^{t} e^{(t-s)\mathcal{L}} \mathcal{PL} e^{s\mathcal{QL}} ds.
    \label{eq:MZ_identity}
\end{equation}
Differentiating \Cref{eq:MZ_identity} with respect to time yields the evolution equation for the full operator:
\begin{equation*}
    \frac{d}{dt} e^{t\mathcal{L}} = e^{t\mathcal{L}} \mathcal{PL} + \int_{0}^{t} e^{(t-s)\mathcal{L}} \mathcal{PL} e^{s\mathcal{QL}} \mathcal{QL} \, ds + e^{t\mathcal{QL}} \mathcal{QL}.
\end{equation*}
Applying this operator equation to the observable $\boldsymbol{\psi}$ and evaluating the result at the initial state $\boldsymbol{x}_0$ yields the Mori-Zwanzig equation, also known as the Generalized Langevin Equation with the distinguished terms:
\begin{equation}
    \frac{\partial}{\partial t} (e^{t\mathcal{L}}\boldsymbol{\psi})(\boldsymbol{x}_0) =
    \underbrace{(e^{t\mathcal{L}} \mathcal{P} \mathcal{L} \boldsymbol{\psi})(\boldsymbol{x}_0)}_{\text{Markovian}} +
    \underbrace{\int_{0}^{t} (e^{(t-s)\mathcal{L}} \mathcal{P} \mathcal{L} e^{s\mathcal{Q}\mathcal{L}} \mathcal{Q} \mathcal{L} \boldsymbol{\psi})(\boldsymbol{x}_0) \, ds}_{\text{Memory}} +
    \underbrace{(e^{t\mathcal{Q}\mathcal{L}} \mathcal{Q} \mathcal{L} \boldsymbol{\psi})(\boldsymbol{x}_0)}_{\text{Noise}}.
    \label{eq:MZ_GLE}
\end{equation}

\Cref{eq:MZ_GLE} that has the following common interpretation:
\begin{enumerate}
    \item \textbf{Markovian term:} the instantaneous part of the dynamics
    projected onto the resolved variables.
    \item \textbf{Memory term:} a convolution integral containing the memory
    kernel $\boldsymbol{K}(t) = \mathcal{P}\mathcal{L} e^{t\mathcal{QL}}
    \mathcal{Q}\mathcal{L}$, which describes how the unresolved dynamics
    influence the current state through their history.
    \item \textbf{Noise term:} the contribution induced
    by the evolution of unresolved initial conditions in the orthogonal
    subspace $\mathcal{Q}$.
\end{enumerate}

\section{Mori-Zwanzig Equation for Linear Systems}\label{sec:linear_MZ}

In this section, we specialize the general Mori-Zwanzig identity to partitioned linear systems. The derivation is organized in three steps: we first define the resolved and unresolved variables, then introduce the projection operator, and finally derive the
exact closed equation for the resolved dynamics. We first define the observed variable assuming this to be the first part of the data vectors. Then we define the projections to model the case, where only the observed part enters the model.

\paragraph{Linear System and State Partition.}
Consider the linear system
\begin{equation}
    \dot{\boldsymbol{x}}(t) = \boldsymbol{A}\boldsymbol{x}(t) + \boldsymbol{b},
    \qquad
    \boldsymbol{x}(0) = \boldsymbol{x}_{0},
    \label{eq:linear_full}
\end{equation}
with $\boldsymbol{x}(t) \in \mathbb{R}^{n}$, $\boldsymbol{A}\in\mathbb{R}^{n\times n}$, and $\boldsymbol{b}\in\mathbb{R}^{n}$. Let $x = (x_1, x_2)$, where $x_1$
is the observed / resolved part.
We split the state and the model into resolved and unresolved components,
\begin{equation}
    \boldsymbol{x} = \begin{bmatrix} \boldsymbol{x}_1 \\ \boldsymbol{x}_2 \end{bmatrix}, \qquad \boldsymbol{A} = \begin{bmatrix} \boldsymbol{A}_{11} & \boldsymbol{A}_{12} \\ \boldsymbol{A}_{21} & \boldsymbol{A}_{22} \end{bmatrix}, \qquad \boldsymbol{b} = \begin{bmatrix} \boldsymbol{b}_{1} \\ \boldsymbol{b}_{2} \end{bmatrix}.
    \label{eq:linear_block_system}
\end{equation}

Our goal is to derive an evolution equation for $\boldsymbol{x}_{1}$ with the help of the Mori-Zwanzig formalism.

\paragraph{Operator Formalism.}
We define the vector of observables simply by the coordinate functions $\boldsymbol{\psi}(\boldsymbol{x}) = \boldsymbol{x}$. The evolution of these observables is governed by the Liouville operator $\mathcal{L}$. For the linear system $\dot{\boldsymbol{x}} = \boldsymbol{A}\boldsymbol{x}+\boldsymbol{b}$, the action of the Liouville operator is
\begin{equation*}
    (\mathcal{L}\psi_i)(\boldsymbol{x})
    = \sum_{k=1}^n (\boldsymbol{A}\boldsymbol{x}+\boldsymbol{b})_k
    \frac{\partial}{\partial x_k}(x_i)
    = \sum_{k=1}^n (\boldsymbol{A}\boldsymbol{x}+\boldsymbol{b})_k \delta_{ik}
    = (\boldsymbol{A}\boldsymbol{x}+\boldsymbol{b})_i.
\end{equation*}
Thus, in vector notation:
\begin{equation*}
    (\mathcal{L}\boldsymbol{\psi})(\boldsymbol{x}) = \boldsymbol{A}\boldsymbol{x}+\boldsymbol{b}.
\end{equation*}

To separate the dynamics, we introduce the projection operator $\mathcal{P}$ and the orthogonal projection $\mathcal{Q} = \mathcal{I} - \mathcal{P}$. The operator $\mathcal{P}$ projects any function of the full state onto the subspace of resolved variables by setting the unresolved variables to zero:
\begin{equation*}
    (\mathcal{P} f)(\boldsymbol{x}) = f(\boldsymbol{x}_1, \boldsymbol{x}_2 = 0).
\end{equation*}
The exact time evolution of the system is described by the Mori-Zwanzig equation:
\begin{equation} \label{eq:MZ_General}
    \frac{\partial}{\partial t}e^{t\mathcal{L}}\boldsymbol{\psi}_{0} =
    \underbrace{e^{t\mathcal{L}}\mathcal{P}\mathcal{L}\boldsymbol{\psi}_{0}}_{\text{Markovian}} +
    \underbrace{e^{t\mathcal{Q}\mathcal{L}}\mathcal{Q}\mathcal{L}\boldsymbol{\psi}_{0}}_{\text{Noise}} +
    \underbrace{\int_{0}^{t}e^{(t-s)\mathcal{L}}\mathcal{P}\mathcal{L}e^{s\mathcal{Q}\mathcal{L}}\mathcal{Q}\mathcal{L}\boldsymbol{\psi}_{0}ds}_{\text{Memory}}.
\end{equation}
Here the subscript $0$ marks the initial time $t=0$, so that $\boldsymbol{\psi}_{0}$ is the observable
at $t=0$ and $e^{t\mathcal{L}}\boldsymbol{\psi}_{0}$ the corresponding observable at time $t$. The same
convention is used for the state components, as in $\boldsymbol{x}_{2,0}=\boldsymbol{x}_{2}(0)$.

\paragraph{Projection Operator and Partial Observation.}

As in the general framework, we choose the coordinate observables
$\boldsymbol{\psi}(\boldsymbol{x}) = \operatorname{id}(\boldsymbol{x}) = \boldsymbol{x}$. For the linear system \Cref{eq:linear_full}, the Liouville operator satisfies
\begin{equation}
    \mathcal{L}(\operatorname{id})(\boldsymbol{x}) = \boldsymbol{A}\boldsymbol{x} + \boldsymbol{b}.
    \label{eq:Lid_linear}
\end{equation}
Assuming that the part $\boldsymbol x_2$ encodes the unresolved coordinates, we
use the projection onto the observable space that simply sets the unresolved
coordinates to zero:
\begin{subequations} \label{eq:linear_projection}
  \begin{align}
    (\mathcal{P}f)(\boldsymbol{x}_{1},\boldsymbol{x}_{2}) &=
    f(\boldsymbol{x}_{1},0), \\
    \intertext{so that for $\mathcal{Q} = \mathcal{I} - \mathcal{P}$, we have}
    (\mathcal{Q}f)(\boldsymbol{x}_{1},\boldsymbol{x}_{2})
    &= f(\boldsymbol{x}_{1},\boldsymbol{x}_{2}) - f(\boldsymbol{x}_{1},0).
  \end{align}
\end{subequations}
 We first apply $\mathcal{L}$ to the identity observable using \Cref{eq:Lid_linear} to obtain
\begin{equation*}
    \mathcal{L}(\operatorname{id})(\boldsymbol{x}) = \boldsymbol{A}\boldsymbol{x} + \boldsymbol{b} =
    \begin{bmatrix}
        \boldsymbol{A}_{11}\boldsymbol{x}_{1} + \boldsymbol{A}_{12}\boldsymbol{x}_{2} + \boldsymbol{b}_1 \\
        \boldsymbol{A}_{21}\boldsymbol{x}_{1} + \boldsymbol{A}_{22}\boldsymbol{x}_{2} + \boldsymbol{b}_2
    \end{bmatrix}.
\end{equation*}
with the partitioned system matrices and vectors as in \Cref{eq:linear_block_system}.
Applying $\mathcal{P}$ and $\mathcal{Q}$ as defined in \Cref{eq:linear_projection},
we obtain
\begin{align*}
    \mathcal{P}\mathcal{L}(\operatorname{id})(\boldsymbol{x}) =
    \begin{bmatrix}
        \boldsymbol{A}_{11}\boldsymbol{x}_{1} + \boldsymbol{b}_1 \\
        \boldsymbol{A}_{21}\boldsymbol{x}_{1} + \boldsymbol{b}_2
    \end{bmatrix}, \qquad
    \mathcal{Q}\mathcal{L}(\operatorname{id})(\boldsymbol{x}) = \mathcal{L}(\operatorname{id})(\boldsymbol{x}) - \mathcal{P}\mathcal{L}(\operatorname{id})(\boldsymbol{x}) =
    \begin{bmatrix}
        \boldsymbol{A}_{12}\boldsymbol{x}_{2} \\
        \boldsymbol{A}_{22}\boldsymbol{x}_{2}
    \end{bmatrix}.
\end{align*}
Here we have introduced the column-split matrices
\begin{equation*}
    \hat{\boldsymbol{A}} =
    \begin{bmatrix}
        \boldsymbol{A}_{11} & \boldsymbol{0} \\
        \boldsymbol{A}_{21} & \boldsymbol{0}
    \end{bmatrix},
    \qquad
    \tilde{\boldsymbol{A}} =
    \begin{bmatrix}
        \boldsymbol{0} & \boldsymbol{A}_{12} \\
        \boldsymbol{0} & \boldsymbol{A}_{22}
    \end{bmatrix},
    \qquad
    \hat{\boldsymbol{A}} + \tilde{\boldsymbol{A}} = \boldsymbol{A},
\end{equation*}
so that the two projected parts take the compact form
\begin{equation*}
    \mathcal{P}\mathcal{L}(\operatorname{id})(\boldsymbol{x})
    = \hat{\boldsymbol{A}}\boldsymbol{x} + \boldsymbol{b},
    \qquad
    \mathcal{Q}\mathcal{L}(\operatorname{id})(\boldsymbol{x})
    = \tilde{\boldsymbol{A}}\boldsymbol{x}.
\end{equation*}
These formulas already show the roles of the two subsystems: $\mathcal{P}\mathcal{L}$
depends only on the resolved state, while $\mathcal{Q}\mathcal{L}$ propagates the
hidden contribution. Note that the constant $\boldsymbol{b}$ appears only in the
projected part, since $\mathcal{Q}$ annihilates constants.


With the projection operator defined, we can construct the Mori-Zwanzig equation
for the linear system with a partially observed and resolved variable, chosen as
$\boldsymbol{x}_1$, by identifying the Markovian, noise, and memory term for the
MZ relation \eqref{eq:MZ_General}.

\paragraph{The Markovian Term.}
The Markovian term represents the instantaneous dynamics projected onto the resolved subspace. From this point on, the observable in \Cref{eq:MZ_General} is the resolved coordinate map
\begin{equation*}
    \boldsymbol{\psi}_{1,0}(\boldsymbol{\xi})
    =
    \boldsymbol{\xi}_{1},
    \qquad
    \boldsymbol{\xi}\in\mathcal{X},
\end{equation*}
where the two subscripts follow the convention already used for $\boldsymbol{x}_{2,0}$: the first
labels the resolved block, as for $\boldsymbol{x}_{1}$, and the second marks the initial time $t=0$.
For this choice, the Markovian contribution is
\begin{equation*}
    \boldsymbol{M}_{1}(t)
    =
    e^{t\mathcal{L}}\mathcal{P}\mathcal{L}\boldsymbol{\psi}_{1,0}.
\end{equation*}
The inner factor is the resolved observable
\begin{align*}
    \boldsymbol{m}_{1}(\boldsymbol{x})
    &:=
    \left(\mathcal{P}\mathcal{L}\boldsymbol{\psi}_{1,0}\right)(\boldsymbol{x})
    =
    \boldsymbol{A}_{11}\boldsymbol{x}_{1}+\boldsymbol{b}_{1}.
\end{align*}
Evaluating the MZ Markovian term at the initial state $\boldsymbol{x}_0$ gives
\begin{equation*}
    \boldsymbol{M}_{1}(t)(\boldsymbol{x}_{0})
    =
    \left(e^{t\mathcal{L}}\boldsymbol{m}_{1}\right)(\boldsymbol{x}_{0})
    =
    \boldsymbol{m}_{1}(\boldsymbol{x}(t))
    =
    \boldsymbol{A}_{11}\boldsymbol{x}_1(t) + \boldsymbol{b}_1.
\end{equation*}

\paragraph{The Noise Term.}
The noise term $\boldsymbol{F}_{1}(t)$ arises from the orthogonal dynamics $\mathcal{Q}\mathcal{L}$. For the resolved observable $\boldsymbol{\psi}_{1,0}(\boldsymbol{\xi})=\boldsymbol{\xi}_{1}$, the contribution in \Cref{eq:MZ_General} is
\begin{equation*}
    \boldsymbol{F}_{1}(t)
    =
    e^{t\mathcal{Q}\mathcal{L}}\mathcal{Q}\mathcal{L}\boldsymbol{\psi}_{1,0}.
\end{equation*}
The inner orthogonal factor is
\begin{align*}
    \boldsymbol{f}_{1,0}(\boldsymbol{x})
    &:=
    \left(\mathcal{Q}\mathcal{L}\boldsymbol{\psi}_{1,0}\right)(\boldsymbol{x})
    =
    \left(\mathcal{L}\boldsymbol{\psi}_{1,0}\right)(\boldsymbol{x})
    -
    \left(\mathcal{P}\mathcal{L}\boldsymbol{\psi}_{1,0}\right)(\boldsymbol{x})
    =
    \boldsymbol{A}_{12}\boldsymbol{x}_{2}.
\end{align*}
The subscripts follow the same convention, so that
$\boldsymbol{f}_{1,0}=\boldsymbol{F}_{1}(0)$ is the resolved orthogonal observable at lag $t=0$.
The resolved constant term $\boldsymbol{b}_{1}$ cancels exactly. Therefore the MZ noise term
evaluated at the initial state is
\begin{equation*}
    \boldsymbol{F}_{1}(t)(\boldsymbol{x}_{0})
    =
    \left(e^{t\mathcal{Q}\mathcal{L}}\boldsymbol{f}_{1,0}\right)(\boldsymbol{x}_{0}).
\end{equation*}
To compute this term explicitly, let
$\boldsymbol{g}(\boldsymbol{x})=\boldsymbol{B}\boldsymbol{x}_{2}$.
Since
\begin{equation*}
    \mathcal{Q}\mathcal{L}\boldsymbol{g}(\boldsymbol{x})
    =
    \boldsymbol{B}\boldsymbol{A}_{22}\boldsymbol{x}_{2},
\end{equation*}
the first powers generated by $\mathcal{Q}\mathcal{L}$ are
\begin{equation*}
    \boldsymbol{f}_{1,0}(\boldsymbol{x}) =
    \boldsymbol{A}_{12}\boldsymbol{x}_{2},
    \quad
    \mathcal{Q}\mathcal{L}\boldsymbol{f}_{1,0}(\boldsymbol{x})
    = \boldsymbol{A}_{12}\boldsymbol{A}_{22}\boldsymbol{x}_{2},
    \quad\text{and}\quad
    (\mathcal{Q}\mathcal{L})^{2}\boldsymbol{f}_{1,0}(\boldsymbol{x})
    =
    \boldsymbol{A}_{12}\boldsymbol{A}_{22}^{2}\boldsymbol{x}_{2}.
\end{equation*}
The pattern gives, by induction,
\begin{equation*}
    (\mathcal{Q}\mathcal{L})^{n}\boldsymbol{f}_{1,0}(\boldsymbol{x})
    =
    \boldsymbol{A}_{12}\boldsymbol{A}_{22}^{n}\boldsymbol{x}_{2},
    \qquad n\ge 0.
\end{equation*}
Thus, through the series expansion of the propagator, we arrive at the explicit
expression for the noise term to read:
\begin{equation*}
    \begin{aligned}
    \boldsymbol{F}_{1}(t)(\boldsymbol{x}_{0})
    = \left(e^{t\mathcal{Q}\mathcal{L}}\boldsymbol{f}_{1,0}\right)(\boldsymbol{x}_{0})
    =
    \left(e^{t\mathcal{Q}\mathcal{L}}
    \mathcal{Q}\mathcal{L}\boldsymbol{\psi}_{1,0}\right)(\boldsymbol{x}_{0})
    = \sum_{n=0}^{\infty}\frac{t^n}{n!}
    \boldsymbol{A}_{12}\boldsymbol{A}_{22}^{n}\boldsymbol{x}_{2,0}
    =
    \boldsymbol{A}_{12} e^{t\boldsymbol{A}_{22}}\boldsymbol{x}_{2,0}.
    \end{aligned}
\end{equation*}

\paragraph{The Memory Term.}
The memory term describes how past states influence the current evolution through a convolution kernel. In \Cref{eq:MZ_General}, this is the contribution
\begin{equation}
    \boldsymbol{G}_{1}(t)
    =
    \int_{0}^{t}
    e^{(t-s)\mathcal{L}}\mathcal{P}\mathcal{L}
    e^{s\mathcal{Q}\mathcal{L}}\mathcal{Q}\mathcal{L}
    \boldsymbol{\psi}_{1,0}\,ds.
    \label{eq:linear_memory_general_term}
\end{equation}
Define the orthogonal-dynamics observable inside the memory integrand by
\begin{equation*}
    \boldsymbol{F}_{1}(s)
    =
    e^{s\mathcal{Q}\mathcal{L}}\mathcal{Q}\mathcal{L}
    \boldsymbol{\psi}_{1,0}.
\end{equation*}
Then the memory term can be rewritten as
\begin{align*}
    \boldsymbol{G}_{1}(t)
    =
    \int_{0}^{t}
    e^{(t-s)\mathcal{L}}\mathcal{P}\mathcal{L}
    \boldsymbol{F}_{1}(s)\,ds,
\end{align*}
that is,
\begin{align*}
    \mathcal{P}\mathcal{L}
    e^{s\mathcal{Q}\mathcal{L}}\mathcal{Q}\mathcal{L}
    \boldsymbol{\psi}_{1,0}
    =
    \mathcal{P}\mathcal{L}\boldsymbol{F}_{1}(s).
\end{align*}
Thus, the main task is to evaluate
$\mathcal{P}\mathcal{L}\boldsymbol{F}_{1}(s)$ and then propagate it with
$e^{(t-s)\mathcal{L}}$. For the resolved observable,
\begin{equation*}
    \boldsymbol{F}_{1}(s)(\boldsymbol{x})
    =
    \left(e^{s\mathcal{Q}\mathcal{L}}\mathcal{Q}\mathcal{L}
    \boldsymbol{\psi}_{1,0}\right)(\boldsymbol{x}).
\end{equation*}
With the closed form of the orthogonal dynamics, treating the noise term at
lag $s$ as an observable of $\boldsymbol{x}$ gives
\begin{equation*}
    \boldsymbol{F}_{1}(s)(\boldsymbol{x})
    =
    \boldsymbol{A}_{12}e^{s\boldsymbol{A}_{22}}\boldsymbol{x}_2.
\end{equation*}
Using the linear form of $\boldsymbol{F}_{1}(s)$, we can write
\begin{equation*}
    \nabla_{\boldsymbol{x}}\boldsymbol{F}_{1}(s)
    =
    \begin{bmatrix}
        \boldsymbol{0} &
        \boldsymbol{A}_{12}e^{s\boldsymbol{A}_{22}}
    \end{bmatrix}.
\end{equation*}
Application of the Liouville operator yields
\begin{align*}
    \mathcal{L}\boldsymbol{F}_{1}(s)(\boldsymbol{x})
    = \nabla_{\boldsymbol{x}}\boldsymbol{F}_{1}(s)(\boldsymbol{x})
    \left(\boldsymbol{A}\boldsymbol{x}+\boldsymbol{b}\right)
    =
    \boldsymbol{A}_{12}e^{s\boldsymbol{A}_{22}}
    \left(
        \boldsymbol{A}_{21}\boldsymbol{x}_1
        + \boldsymbol{A}_{22}\boldsymbol{x}_2
        + \boldsymbol{b}_2
    \right).
\end{align*}
After projection,
\begin{align*}
    \mathcal{P}\mathcal{L}\boldsymbol{F}_{1}(s)(\boldsymbol{x})
    = \left(\mathcal{L}\boldsymbol{F}_{1}(s)\right)(\boldsymbol{x}_1,\boldsymbol{0})
    =
    \boldsymbol{A}_{12}e^{s\boldsymbol{A}_{22}}
    \left(
        \boldsymbol{A}_{21}\boldsymbol{x}_1
        + \boldsymbol{A}_{22}\boldsymbol{0}
        + \boldsymbol{b}_2
    \right)=
    \boldsymbol{A}_{12}e^{s\boldsymbol{A}_{22}}
    \left(
        \boldsymbol{A}_{21}\boldsymbol{x}_1 + \boldsymbol{b}_2
    \right).
\end{align*}
Substitution of this expression back into the integrand in
\Cref{eq:linear_memory_general_term} gives
\begin{align*}
    \left(
    e^{(t-s)\mathcal{L}}\mathcal{P}\mathcal{L}
    e^{s\mathcal{Q}\mathcal{L}}\mathcal{Q}\mathcal{L}
    \boldsymbol{\psi}_{1,0}
    \right)(\boldsymbol{x}_{0})
    =
    \left(
    e^{(t-s)\mathcal{L}}\mathcal{P}\mathcal{L}
    \boldsymbol{F}_{1}(s)
    \right)(\boldsymbol{x}_{0})=
    \boldsymbol{A}_{12}e^{s\boldsymbol{A}_{22}}
    \left(
        \boldsymbol{A}_{21}\boldsymbol{x}_1(t-s)
        + \boldsymbol{b}_2
    \right).
\end{align*}
Therefore,
\begin{align*}
    \boldsymbol{G}_1(t)
    =
    \int_0^t
    \boldsymbol{A}_{12}e^{s\boldsymbol{A}_{22}}
    \left(
        \boldsymbol{A}_{21}\boldsymbol{x}_1(t-s)
        + \boldsymbol{b}_2
    \right) ds
    =
    \int_0^t
    \boldsymbol{A}_{12}e^{s\boldsymbol{A}_{22}}
    \boldsymbol{A}_{21}\boldsymbol{x}_1(t-s)\,ds
    +
    \int_0^t
    \boldsymbol{A}_{12}e^{s\boldsymbol{A}_{22}}
    \boldsymbol{b}_2\,ds.
\end{align*}
The first integral is the history-dependent part of the memory: the current resolved
rate depends on the past values $\boldsymbol{x}_1(t-s)$ through a coefficient
that depends only on the elapsed time $s$. Defining the memory kernel by
\begin{equation*}
    \boldsymbol{K}(s) = \boldsymbol{A}_{12} e^{s\boldsymbol{A}_{22}} \boldsymbol{A}_{21}, \qquad s \ge 0,
\end{equation*}
this history-dependent part can be written compactly as
\begin{equation*}
    \int_0^t \boldsymbol{K}(s)\boldsymbol{x}_1(t-s)\,ds.
\end{equation*}
The second integral has the same unresolved propagator but contains no history
of $\boldsymbol{x}_1$. In the final MZ form below, we keep this contribution
inside the memory term, since it is still generated through the unresolved
block.


Combining the Markovian, noise, and memory terms derived from
$\boldsymbol{\psi}_{1,0}(\boldsymbol{\xi})=\boldsymbol{\xi}_{1}$ yields
\begin{align}
    \dot{\boldsymbol{x}}_1(t) ={}& \boldsymbol{A}_{11}\boldsymbol{x}_1(t)
    + \boldsymbol{b}_{1}
    + \boldsymbol{A}_{12} e^{t\boldsymbol{A}_{22}} \boldsymbol{x}_{2}(0)
    + \int_{0}^{t} \boldsymbol{A}_{12} e^{s\boldsymbol{A}_{22}}
    \left(
        \boldsymbol{A}_{21}\boldsymbol{x}_{1}(t-s)+\boldsymbol{b}_{2}
    \right) ds,
    \label{eq:linear_MZ_final}
\end{align}
\Cref{eq:linear_MZ_final} is exactly the Mori-Zwanzig equation for the resolved observable of the linear system.


Comparing \Cref{eq:linear_MZ_final} with the general MZ  \Cref{eq:MZ_GLE}, we identify the Markovian, noise, and memory contributions for the resolved observed state $\boldsymbol{x}_1$ equation:
\begin{align}
    (\text{Markovian})_1 &= \boldsymbol{A}_{11}\boldsymbol{x}_{1}(t) + \boldsymbol{b}_1,
    \notag \\
    (\text{Noise})_1 &= \boldsymbol{A}_{12} e^{t\boldsymbol{A}_{22}} \boldsymbol{x}_{2}(0),
    \notag \\
    (\text{Memory})_1 &=
    \int_{0}^{t}
    \boldsymbol{A}_{12} e^{s\boldsymbol{A}_{22}}
    \left(
        \boldsymbol{A}_{21}\boldsymbol{x}_{1}(t-s)+\boldsymbol{b}_{2}
    \right) ds
    =
    \int_{0}^{t} \boldsymbol{K}(s)\boldsymbol{x}_{1}(t-s)\, ds
    +
    \int_{0}^{t} \boldsymbol{A}_{12} e^{s\boldsymbol{A}_{22}}
    \boldsymbol{b}_{2}\, ds,
    \label{eq:linear_memory}
\end{align}
with memory kernel
\begin{equation}
    \boldsymbol{K}(s) =
    \boldsymbol{A}_{12} e^{s\boldsymbol{A}_{22}} \boldsymbol{A}_{21},
    \qquad s \ge 0.
    \label{eq:linear_memory_kernel}
\end{equation}
This $\boldsymbol{x}_1$ decomposition makes the interpretation transparent: $\boldsymbol{A}_{11}\boldsymbol{x}_1+\boldsymbol{b}_1$ governs the instantaneous resolved dynamics, the term
$\boldsymbol{A}_{12} e^{t\boldsymbol{A}_{22}} \boldsymbol{x}_{2}(0)$ carries the effect of the hidden initial condition, and the memory term contains both the history feedback through the kernel $\boldsymbol{K}$ and the forcing from the hidden block.

\begin{remark}
We note that the same formulas can also be derived with the help of the
variation-of-constants formula for the unresolved component. However, the
Mori-Zwanzig derivation via observables and projectors is the one that extends
naturally to nonlinear systems. To see that, we point out that all the derivations
of the formulas in \Cref{sec:MZ_formalism} are for general nonlinear systems.
\end{remark}

\begin{remark}
We note that \Cref{eq:linear_MZ_final} is exact but does not yet give a model for the observed
variables alone: namely the memory kernel and noise term still involve the operators of the
unobserved dynamics, and the noise term still depends on the hidden
initial condition. These are two separate issues that need to be addressed in
practical applications.
Firstly, we think that the operators are not so much a concern as the
kernels of the reduced model need to be inferred or approximated anyways.
Secondly, the hidden initial
condition only enters through the noise term so that it can treated together
with possible perturbation. Alternatively, it may
be estimated from the observed trajectory data with classical techniques from
observer design.
\end{remark}

\section{Illustrative Example: Mori-Zwanzig for Partial Observation of the Harmonic Oscillator}

We consider a classical harmonic oscillator with natural frequency $\omega$ and
the
state at time $t$ described by the vector $\boldsymbol{u}(t) = (x(t), v(t))$, consisting of
the position $x(t)$ and velocity $v(t)$, and the equations of motion given
through the linear time-invariant system

\begin{equation*}
    \dot{\boldsymbol{u}} =
    \begin{bmatrix}
        0 & 1 \\
        -\omega^2 & 0
    \end{bmatrix}
    \boldsymbol{u}.
\end{equation*}

We now apply the general linear MZ result of
\Cref{eq:linear_memory,eq:linear_memory_kernel} directly to this example. We
identify the resolved and hidden variables and the projected parts of the
systems as
\begin{equation*}
    \boldsymbol{x}_1 = x,
    \quad
    \boldsymbol{x}_2 = v,
    \quad\text{and }
    \boldsymbol{A}_{11}=0,\quad
    \boldsymbol{A}_{12}=1,\quad
    \boldsymbol{A}_{21}=-\omega^2,\quad
    \boldsymbol{A}_{22}=0,
    \quad
    \boldsymbol{b}_1=\boldsymbol{b}_2=0,
\end{equation*}
respectively.
Substituting this into the general formulas, we obtain the three MZ terms directly. The
Markovian term vanishes since $\boldsymbol{A}_{11}=0$ and $b_1=0$:
\begin{equation}
    (\text{Markovian})_x
    =
    \boldsymbol{A}_{11}x(t)+\boldsymbol{b}_1
    =
    0,
    \label{eq:markovian}
\end{equation}
whereas the noise term carries the hidden initial condition $v_0=\boldsymbol{x}_2(0)$:
\begin{equation}
    (\text{Noise})_x
    =
    \boldsymbol{A}_{12}\,e^{t\boldsymbol{A}_{22}}\boldsymbol{x}_2(0)
    = 1 e^0 v_0 =
    v_0.
    \label{eq:noise}
\end{equation}
Finally, with the memory kernel \Cref{eq:linear_memory_kernel} as
\begin{equation*}
    \boldsymbol{K}(s)
    =
    \boldsymbol{A}_{12}\,e^{s\boldsymbol{A}_{22}}\boldsymbol{A}_{21}
    =
    -\omega^2,
\end{equation*}
the memory term becomes
\begin{equation}
    (\text{Memory})_x
    =
    \int_0^t \boldsymbol{K}(s)\,x(t-s)\,ds
    =
    -\omega^2\int_0^t x(t-s)\,ds.
    \label{eq:memory}
\end{equation}

Combining \Cref{eq:markovian,eq:noise,eq:memory}, we obtain the Mori-Zwanzig equation for the resolved variable $x$
\begin{equation}
    \dot{x}(t)
    =
    (\text{Markovian})_x + (\text{Noise})_x + (\text{Memory})_x
    =
    0 + v_0 - \omega^2\int_0^t x(t-s)\,ds.
    \label{eq:harmonic_mz_reduced}
\end{equation}


\begin{remark}
  We briefly illustrate the consistency of the MZ approach with an ad-hoc
  derivation of the governing equation for the position of the mass: If we
  differentiate \Cref{eq:harmonic_mz_reduced} with respect to time, we
  immediately arrive at
  the standard formulation as second order system
\begin{equation*}
    \ddot{x}(t) = - \omega^2 x(t),
\end{equation*}
which, with initial data $x(0)=x_0$ and $\dot{x}(0)=v_0$, has the well-known
explicit solution as
\begin{equation*}
    x(t) = x_0\cos(\omega t) + \frac{v_0}{\omega}\sin(\omega t).
\end{equation*}
\end{remark}


\section{Numerical Examples}\label{sec:Numerical_example}
We validate the MZ formalism on two examples: a damped oscillator chain and a 2D wave equation, each with partial observation.
In the numerical realization, solving the MZ equation requires the discretization of the resulting integro-differential equation.
We use the trapezoidal rule to approximate the memory integral.
The reference trajectories are computed with the \texttt{solve\_ivp} routine of the SciPy library \cite{SciPy2020} using RK45, a high-order Runge--Kutta integration scheme,
and are used as the reference solution for the present comparison.
Since the MZ equation is exact at the continuous level, the observed differences mainly reflect the numerical error
introduced by the time discretization and by the numerical approximation of the memory integral.
As expected, these errors are of the order of the chosen time-step size.
Therefore, we focus below on the simulation results and do not provide a separate error analysis.

\paragraph{Damped Oscillatory System.}
We consider a chain of $N_s = 5$ oscillators with uniform mass $m = 1.0$ kg, spring constant $k = 1.0$ N/m, and damping coefficient $c = 0.1$ N\,s/m. The system is discretized in time using $\Delta t = 0.01$ s. External periodic forcing $f_i(t) = F_0 \sin(\omega_f t + \phi_i)$ with $F_0 = 0.5$ N and $\omega_f = 1.5$ rad/s is applied to selected oscillators, where $\phi_i$ is the prescribed phase offset for oscillator $i$. This forcing is treated as a prescribed time-dependent input $\boldsymbol{b}(t)$ in the linear system description.

We observe $m_s = 3$ oscillators ($m=6$ states: positions and velocities), while $N_s - m_s = 2$ oscillators remain unobserved. The simulation spans $T = 500$ time steps ($t_{\text{final}} = 5.0$ s). Reference trajectories are generated using \texttt{solve\_ivp} with RK45.

\begin{figure}
    \begin{center}
        \includegraphics[width=0.9\textwidth]{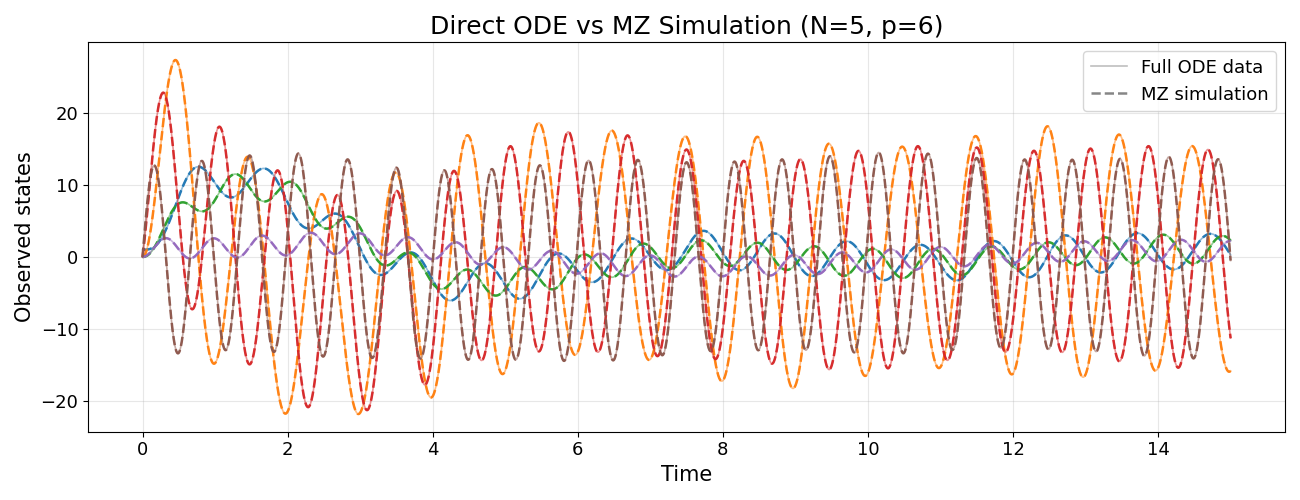}
        \caption{Trajectories of the reference solution (solid) and MZ solution (dashed) for the six observed state variables of the damped oscillator chain under partial observation ($m_s=3$ of $N_s=5$ oscillators observed).} \label{fig:damped_ODE}
    \end{center}
\end{figure}

The comparison between the MZ solution with partial observation and the full reference solution is shown in \Cref{fig:damped_ODE}. The plot shows the trajectories of the reference solution and the MZ solution for six state variables. The curves agree up to the integration error introduced by the discretized MZ equation.

\paragraph{2D Wave Equation.}
We consider a damped 2D wave equation on $\Omega = [0,L_x] \times [0,L_y]$ with homogeneous Dirichlet boundary conditions:
\begin{equation*}
    \frac{\partial^2 u}{\partial t^2} + \gamma \frac{\partial u}{\partial t}
    = c^2 \left( \frac{\partial^2 u}{\partial x^2} + \frac{\partial^2 u}{\partial y^2} \right),
    \quad (x,y) \in \Omega.
\end{equation*}
Introducing $v = \partial_t u$, the first-order form is
\begin{equation*}
    \frac{d}{dt}
    \begin{bmatrix}
        \boldsymbol{u} \\
        \boldsymbol{v}
    \end{bmatrix}
    =
    \begin{bmatrix}
        \boldsymbol{0} & \boldsymbol{I} \\
        c^2\boldsymbol{L} & -\gamma \boldsymbol{I}
    \end{bmatrix}
    \begin{bmatrix}
        \boldsymbol{u} \\
        \boldsymbol{v}
    \end{bmatrix},
\end{equation*}
where $\boldsymbol{L} = \boldsymbol{I}_y \otimes \boldsymbol{T}_x + \boldsymbol{T}_y \otimes \boldsymbol{I}_x$ is the discrete Laplacian assembled via second-order central differences (five-point stencil). The spatial domain is discretized on a uniform $n_x \times n_y = 31 \times 31$ grid with $L_x = L_y = 1.0$, giving $n = 961$ spatial unknowns and state dimension $2n = 1922$. We set $c=0.5$, $\gamma=0.5$, and simulate until $T=1.0$ with $N_t=10{,}000$ time steps ($\Delta t \approx 10^{-4}$) using \texttt{solve\_ivp} with RK45. The initial displacement is a centered Gaussian pulse,
\begin{equation*}
    u_0(x,y)=\exp\left(-100\left[(x-0.5)^2+(y-0.5)^2\right]\right),
\end{equation*}
with zero initial velocity. We observe three quadrants of the spatial domain, leaving one quadrant unobserved.

\begin{figure}[htbp]
    \centering
        \includegraphics[width=0.74\textwidth]{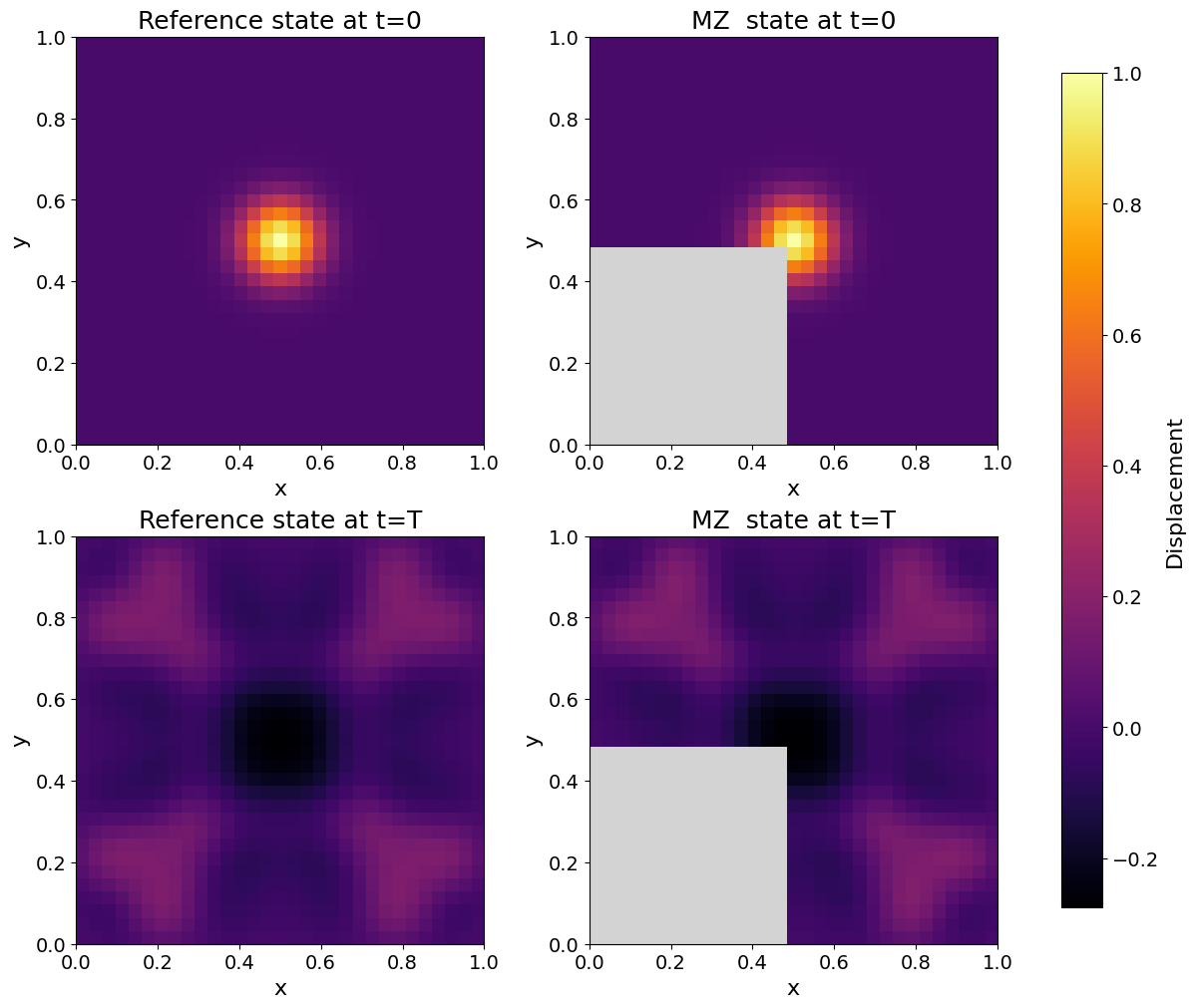}
        \caption{The results of the 2D wave equation modelled by MZ with some
        parts (in gray) being unobserved. The left panel illustrates the reference solution, while the right
panel demonstrates the future prediction of the observed variables using the
reconstructed MZ model. The upper row shows the initial state in the simulation,
and the lower row shows the state at final time.
}\label{fig:2d_wave_31}
\end{figure}

As shown in \Cref{fig:2d_wave_31}, the MZ reconstruction and future prediction
of the observed variables for the 2D wave equation closely match the reference
solution.

\section{Conclusions}\label{sec:conclusion}
This paper presents a comprehensive treatment of the Mori-Zwanzig formalism for linear dynamical systems under partial observation. We derived exact closed-form expressions for the memory kernel, Markovian term, noise term, and forcing contribution, providing analytical clarity that is typically unavailable in nonlinear settings. The harmonic oscillator example illustrated how projection operators systematically redistribute hidden variable effects into memory integrals and initial-condition terms, maintaining exact equivalence with the full dynamics.

Numerical experiments on a damped oscillator chain and a 2D wave equation
demonstrate how the proposed formulation reproduces the reference trajectories
for observed variables despite observing only a subset of state variables or
spatial regions. The remaining discrepancies are governed by the discretization
of the MZ integro-differential equation and by the trapezoidal approximation of
the memory integral, while the continuous MZ formulation itself is exact.\nopagebreak[4]
The linear MZ framework established is supposed to serve as a rigorous benchmark for
extending and understanding the methodology for nonlinear systems. Future work will investigate
also approximative nonlinear extensions of the formalism, spectral and data-driven
approximations of the memory kernel, and efficient numerical
solution of the resulting reduced dynamical equations.


\section*{Acknowledgments}%
\addcontentsline{toc}{section}{Acknowledgments}
  Fan Wang gratefully acknowledges the German Research Foundation (DFG) through the research training group 2297 ``MathCoRe'', Magdeburg.


\addcontentsline{toc}{section}{References}
\bibliographystyle{plain}
\bibliography{reference}

\end{document}